\documentclass[journal,twoside,web]{ieeecolor}

\usepackage{cite}
\usepackage{amsmath,amssymb,amsfonts}
\usepackage{algorithmic}
\usepackage{graphicx}
\usepackage{textcomp}

\usepackage{hyperref}
\hypersetup{hidelinks}
\usepackage{mathtools}
\mathtoolsset{showonlyrefs}

\usepackage{tikz}
\usepackage{pgfplots}
\pgfplotsset{compat=1.18}
\usepgfplotslibrary{groupplots}

\usepackage{graphicx}
\usepackage{subfigure}
\usepackage{epsfig}
\usepackage{cite}
\usepackage{lcsys}

\newtheorem{theorem}{Theorem}[section]
\newtheorem{lemma}[theorem]{Lemma}
\newtheorem{proposition}[theorem]{Proposition}

\newtheorem{definition}[theorem]{Definition}

\DeclareMathOperator{\dom}{dom}
\DeclareMathOperator{\interior}{int}
\DeclareMathOperator{\arcsinh}{arcsinh}

\DeclareMathOperator{\sech}{sech}

\DeclareMathOperator{\blkdiag}{blkdiag}
\DeclareMathOperator*{\minimize}{minimize}
\DeclareMathOperator*{\st}{subject\,to}

\allowdisplaybreaks
\begin{document}

\def\BibTeX{{\rm B\kern-.05em{\sc i\kern-.025em b}\kern-.08em
    T\kern-.1667em\lower.7ex\hbox{E}\kern-.125emX}}
\markboth{\journalname, VOL. XX, NO. XX, XXXX 2017}
{Author \MakeLowercase{\textit{et al.}}: Preparation of Papers for IEEE Control Systems Letters (August 2022)}

\title{Mirror-Space Constrained Optimal Control}

\author{
Casian Iacob$^{1, 2}$,
Hany Abdulsamad$^{3}$,
Simo S\"arkk\"a$^{1, 2}$, \IEEEmembership{Senior Member, IEEE}
\thanks{$^{1}$Department of Electrical Engineering and Automation, Aalto University, 02150 Espoo, Finland; \texttt{\{casian.iacob, simo.sarkka\}@aalto.fi}}
\thanks{$^{2}$ELLIS Institute Finland}
\thanks{$^{3}$Amsterdam Machine Learning Lab, University of Amsterdam, 1098 XH Amsterdam, The Netherlands; \texttt{h.abdulsamad@uva.nl}}
}

\maketitle
\thispagestyle{empty}

\begin{abstract}
We present a method for nonlinear optimal control problems with convex control constraints. In the proposed approach, we use the gradient of a Legendre-type function to map the controls to an unconstrained mirror space, where the problem is solved with the iterative linear quadratic regulator. Using convex conjugacy, we express all derivatives explicitly and avoid differentiation through the inverse map. Feasibility holds by construction without additional hyperparameters. We demonstrate in numerical experiments that this method requires fewer iterations than interior-point methods.
\end{abstract}


\section{Introduction}
\label{sec:introduction}

\IEEEPARstart{O}{ptimal control} problems~\cite{lewis2012optimal} subject to inequality constraints underpin many applications in  robotics and control systems design~\cite{rawlings2017model, lynch2017modern}. Model predictive control (MPC)~\cite{rawlings2017model} handles such constraints by repeatedly solving a finite-horizon constrained optimal control problem online. In this paper, we address a specific class of optimal control problems subject to convex control constraints for MPC.

Various strategies have been proposed to incorporate constraint handling into the dynamic programming framework \cite{tassa2012synthesis} for MPC. For control-limited problems, the method in~\cite{tassa2014control} solves a quadratic program at each stage of the backward pass. More general constraints are handled via penalty or barrier formulations~\cite{howell2019altro, ma2022alternating, chen2019autonomous}, primal-dual interior-point (IP) methods~\cite{hansson2000primal, pavlov2021interior}, or  augmented Lagrangian (AL) methods~\cite{gill2012primal, jallet2022constrained, ma2022alternating, rey2020admm, howell2019altro}. These approaches often require  heuristic barrier or penalty parameter tuning, and for IP methods, ensuring the positivity of slack variables, for example through the fraction-to-the-boundary rule~\cite{nocedal2006numerical}. Furthermore, AL methods~\cite{ma2022alternating, rey2020admm, howell2019altro}, can suffer from slow convergence. Recent work has proposed mapping state and input variables to an unconstrained space via diffeomorphisms~\cite{zhu2025diffeomorphism}. However, the resulting backward pass requires differentiating through the inverse maps, which can become intractable for coupled convex sets.

Concurrently, there is growing interest in generalizing optimization methods to non-Euclidean geometries by replacing Euclidean metrics with Bregman divergences~\cite{wang2014bregman, hajar2025beyond, hassibi2025beyond}. Related work develops variational formulations of mirror descent and predictive mirror dynamics~\cite{tzen2023variational, pan2026unified}, and connects MPC to online learning through dynamic mirror descent~\cite{wagener2019online}. These approaches exploit mirror geometry and Bregman divergences to characterize dynamics or structure optimization updates. In contrast, we use the mirror map to reparameterize the control variables of the optimal control problem itself.

The contribution of this paper is an iterative dynamic programming algorithm, based on mirror maps, for nonlinear optimal control problems with convex control constraints. We transform such problems into unconstrained problems in a mirror space through Legendre-type distance-generating functions (DGFs) and their associated mirror maps~\cite{bubeck2015convex}. For log-barrier DGFs, the mirror map is closely related to the barrier functions of interior-point methods. However, its effect is realized through a change of variables instead of augmenting the cost with a weighted penalty term, which traditionally requires a heuristic annealing schedule~\cite{nocedal2006numerical, pavlov2021interior}. We adapt the unconstrained iLQR \cite{tassa2012synthesis} to solve the transformed problem. Using convex conjugacy, we express the derivatives required in the backward pass in terms of the control-space variables, avoiding differentiation through potentially intractable inverse maps. Unlike penalty, barrier, and augmented Lagrangian approaches, the proposed method introduces no hyperparameters beyond those of standard iLQR~\cite{tassa2012synthesis}.

\section{Problem Formulation}
\label{sec:problem-formulation}

We consider the finite-horizon optimal control problem
\begin{equation}
    \label{eq:optimal-control-problem}
    \begin{split}
        \minimize_{u_{0:T-1}} \quad &\mathcal{J}(u_{0:T-1}) = c_T(x_T) + \sum_{t=0}^{T-1} c_t(x_t, u_t),\\
        \st \quad & x_{t+1} = f(x_t, u_t), \quad x_0 = x_s, \\
                  & u_t \in \mathcal{U}, \quad t=0,\ldots,T-1,
    \end{split}
\end{equation}
where $x_t\in\mathbb{R}^{n}$ and $u_t\in\mathbb{R}^{m}$ are the state and control, respectively. The initial state is $x_s\in\mathbb{R}^{n}$ and $T$ is the horizon. The admissible control set is $\mathcal{U} = \{u\in\mathbb{R}^{m}\colon g(u)\leq 0\}$, where $\mathcal{U}$ is closed, convex, and has nonempty interior. The vector-valued function $g$ is continuous and componentwise convex, the dynamics $f$ and the costs $c_t$ and $c_T$ are twice continuously differentiable, and $\mathcal{J}$ is bounded from below.

In the following section, we introduce a new dynamic programming approach that handles constraints by reparameterizing the optimal control problem through mirror maps and applying unconstrained iLQR in the resulting mirror space.

\begin{table*}[t]
    \centering
    \caption{Examples of Mirror Maps and Corresponding Constraints. {\protect\unboldmath Here $\mu=(b+a)/2$ and $\rho=(b-a)/2$.}}
    \label{tab:mm-table}
    \normalsize
    \renewcommand{\arraystretch}{1.4}
    \begin{tabular}{lcccc}
        \hline
        Name & $\interior(\mathcal{U})$ & DGF $\psi(u)$ & Map $\nabla \psi (u)$ & Inverse $\nabla \psi^*(y)$\\
        \hline
        Box (entropy) & $a < u < b$ & $(b-u)\log(b-u) + (u-a)\log(u-a)$ & $\log\frac{u-a}{b-u}$ & $\mu+\rho\tanh\frac{y}{2}$\\
        Box (log-barrier) & $a < u < b$ & $-\log(b-u) - \log(u-a)$ & $\frac{1}{b-u} - \frac{1}{u-a}$ & $\mu+\rho\tanh\big(\frac{1}{2}\arcsinh(\rho y)\big)$\\
        Positivity & $u > 0$ & $\sum_i (u^i\log u^i - u^i)$ & $\log u$ & $e^{y}$\\
        Euclidean ball & $\|u\|_2 < 1$ & $-\sqrt{1-\|u\|_2^2}$ & $\frac{u}{\sqrt{1-\|u\|_2^2}}$ & $\frac{y}{\sqrt{1+\|y\|_2^2}}$\\
        \hline
    \end{tabular}
\end{table*}

\section{Mirror-Space Optimal Control}
\label{sec:mmilqr}

In this section, we map the constrained optimal control problem to an unconstrained mirror space and derive the dynamic programming recursion directly in mirror coordinates.

\subsection{Mirror Maps for Convex Constraints}

To formalize the mirror maps used in our construction, we recall the notion of Legendre-type functions.

\begin{definition}[Legendre-type function]
    A proper, closed, convex function $\psi\colon\mathbb{R}^m\to(-\infty,\infty]$ whose domain has nonempty interior, $\interior(\dom\psi)\neq\emptyset$, is of \emph{Legendre type} if it is twice continuously differentiable on $\interior(\dom\psi)$ with positive definite Hessian $\nabla^2\psi(u)\succ 0$ for all $u\in\interior(\dom\psi)$, and $\|\nabla\psi(u_k)\|\to\infty$ for every sequence $(u_k)\subset\interior(\dom\psi)$ converging to a boundary point of $\dom\psi$.
\end{definition}

The following proposition relates the gradients and Hessians of a Legendre-type function and its convex conjugate.

\begin{proposition}[Global diffeomorphism]
    \label{prop:global-diffeo}
    Let $\psi$ be of Legendre type and let $\psi^*$ be its convex conjugate, $\psi^*(y) = \sup_u\,\left\{ y^\top u-\psi(u)\right\}$. Then, $\nabla\psi\colon\interior(\dom\psi) \to \interior(\dom\psi^*)$ is a diffeomorphism with inverse $\nabla \psi^*$, and $\nabla^2 \psi^*(y) = \nabla^2 \psi(u)^{-1}$, $\forall u\in\interior(\dom\psi)$, and $y = \nabla\psi(u)$.
\end{proposition}

\begin{proof}
    By~\cite[Thrm.~26.5]{rockafellar1997convex}, $\nabla\psi$ is a bijection from $\interior(\dom\psi)$ onto $\interior(\dom\psi^*)$ with inverse $\nabla\psi^*$. Since $\nabla\psi$ is continuously differentiable with nonsingular Jacobian $\nabla^2\psi(u)\succ 0$, the inverse function theorem shows that $\nabla\psi^*$ is continuously differentiable with Jacobian $\nabla^2\psi^*(y) = \partial u/\partial y = (\partial y/\partial u)^{-1} = \nabla^2\psi(u)^{-1}$.
\end{proof}

We choose a Legendre-type DGF $\psi$ with $\interior(\dom\psi)=\interior(\mathcal{U})$ and $\dom\psi^*=\mathbb{R}^{m}$. The latter is automatic for a bounded $\mathcal{U}$ and can be checked directly otherwise. Proposition~\ref{prop:global-diffeo} therefore gives $\nabla\psi\colon \interior(\mathcal U) \to \mathbb{R}^{m}$ and $\nabla\psi^*\colon \mathbb{R}^{m} \to \interior(\mathcal U)$. The function $\psi$ acts as a potential over the domain. It is smooth on the interior and steep at the boundary, ensuring that $\nabla\psi^*$ maps every $y\in\mathbb{R}^m$ into $\interior(\mathcal{U})$.

At each stage $t$, the mirror map is constructed to transform the constrained control $u_t$ into an unconstrained mirror control $y_t = \nabla \psi (u_t)$, where $y_t \in \mathbb{R}^{m}$. By Proposition~\ref{prop:global-diffeo}, the inverse $\nabla \psi^*$ maps $y_t$ back to a control in the interior of $\mathcal U$, $u_t = \nabla \psi^*(y_t)$. Feasibility therefore holds by construction.

Table~\ref{tab:mm-table} lists mirror maps with tractable inverse maps. Similar constructions approximate constraint indicator functions in~\cite{lutter2020hjb}.  Here, they are used instead as coordinate transformations to an unconstrained mirror space.

\subsection{Dynamic Programming in Mirror Coordinates}
\label{sec:mirror-dp}
Substituting $u_t=\nabla\psi^*(y_t)$ in problem~\eqref{eq:optimal-control-problem} yields
\begin{equation}
    \label{eq:mirror-ocp}
    \begin{split}
        \minimize_{y_{0:T-1}} \quad & c_T(x_T) + \sum_{t=0}^{T-1} c_t\big(x_t,\nabla\psi^*(y_t)\big) \\
        \st \quad & x_{t+1}=f\big(x_t,\nabla\psi^*(y_t)\big), \quad  x_0 = x_s.
    \end{split}
\end{equation}
Since $\nabla\psi^*$ maps $\mathbb{R}^m$ onto $\interior(\mathcal{U})$, problem~\eqref{eq:mirror-ocp} is equivalent to~\eqref{eq:optimal-control-problem} with $\mathcal{U}$ replaced by $\interior(\mathcal{U})$, and all iterates are strictly feasible. As $\mathcal{U}$ is convex with nonempty interior, it is the closure of $\interior(\mathcal{U})$, so by continuity of $\mathcal{J}$ both problems have the same optimal value. If~\eqref{eq:optimal-control-problem} has a minimizer with all controls in $\interior(\mathcal{U})$, it is attained in the mirror space. If constraints are active at every minimizer, the optimal value is approached as $\|y_t\|\to\infty$, when the controls approach the boundary from the interior, similar to interior-point methods.

The mirror controls $y_t$ are unconstrained, allowing us to apply dynamic programming to problem~\eqref{eq:mirror-ocp}. Starting from $V_T(x_T) = c_T(x_T)$, the Bellman recursion~\cite{lewis2012optimal} reads
\begin{align}
    Q_t(x_t,u_t) &\coloneqq c_t(x_t,u_t) + V_{t+1}\big(f(x_t,u_t)\big),\label{eq:q-function}\\
    V_t(x_t) &=\min_{y_t} \, Q_t\big(x_t,\nabla\psi^*(y_t)\big),\label{eq:bellman}
\end{align}
for $t=0,\ldots,T-1$.
The mirror map enters only through the minimization in~\eqref{eq:bellman}, which is carried out over $y_t\in\mathbb{R}^m$ and thus implicitly over $u_t\in\interior(\mathcal{U})$.

Solving~\eqref{eq:bellman} requires derivatives of functions of the form $h(y) = g(\nabla\psi^*(y))$ with respect to the mirror control. Computing them directly requires the first and second derivatives of the inverse map $\nabla\psi^*$, which may not be available in closed form. Instead, we express them in terms of the derivatives of $\psi$ only. To this end, let $\Psi(u) = \nabla^2\psi(u)$ denote the Hessian of the DGF, assume that $\psi\in C^3(\interior(\mathcal{U}))$, and let $\mathcal{T}_\psi(u)$ denote the third-order derivative tensor with elements
\begin{equation}
    [\mathcal{T}_\psi(u)]_{ijk} = \frac{\partial^3 \psi(u)}{\partial u^i\partial u^j \partial u^k}.\label{eq:third-order-tensor}
\end{equation}
Contracting $\mathcal{T}_\psi(u)$ with a vector $v$ yields a matrix, which we write as $\mathcal{T}_\psi(u) \cdot v$ with elements
\begin{equation}
    [\mathcal{T}_\psi(u) \cdot v]_{ij} = \sum_k\frac{\partial^3\psi(u)}{\partial u^i\partial u^j\partial u^k}v^k.\label{eq:tensor-dot}
\end{equation}

\begin{lemma}[Derivatives under the mirror map]
\label{lem:mirror-derivatives}
Let $g$ be a twice continuously differentiable function on $\interior(\mathcal{U})$, and let $h(y) = g(\nabla\psi^*(y))$. Then, for every $y\in\mathbb{R}^m$ and $u=\nabla\psi^*(y)$,
\begin{align}
    \nabla h(y) &= \Psi(u)^{-1}\,\nabla g(u),\label{eq:mirror-gradient}\\
    \nabla^2 h(y) &= \Psi(u)^{-1}\big[\nabla^2 g(u) - \mathcal{T}_\psi(u)\cdot\nabla h(y)\big]\Psi(u)^{-1}.\label{eq:mirror-hessian}
\end{align}
\end{lemma}
\begin{proof}
    See Appendix~\ref{app:proof-mirror-derivatives}.
\end{proof}
Lemma~\ref{lem:mirror-derivatives} involves only the Hessian and the third-order derivative of $\psi$, evaluated at $u = \nabla\psi^*(y)$.

\subsection{Mirror-Space iLQR}
\label{sec:mirror-ilqr}
In this section, we apply the iLQR approximation~\cite{tassa2012synthesis} to the Bellman recursion~\eqref{eq:bellman} and relate the resulting mirror-space iteration to iLQR in control space. We expand the $Q$-function and the value function to second order around a nominal trajectory $\bar{x}_{0:T}$, $\bar{u}_{0:T-1}$ with $\bar{u}_t\in\interior(\mathcal{U})$ and $\bar{y}_t = \nabla\psi(\bar{u}_t)$:
\begin{equation}
    \begin{aligned}
        \delta Q_t(x_t, u_t)
        &\approx
          \delta x_t^\top Q_{x,t}
          + \delta u_t^\top Q_{u,t}
          + \delta x_t^\top Q_{xu,t}\delta u_t\\
        &\quad +\tfrac{1}{2}\delta x_t^\top Q_{xx,t}\delta x_t + \tfrac{1}{2}\delta u_t^\top Q_{uu,t}\delta u_t, \\
        \delta V_t(x_t)
        &\approx  \delta x_t^\top V_{x,t} + \tfrac{1}{2}\delta x_t^\top V_{xx,t}\delta x_t,
    \end{aligned}
\end{equation}
where $\delta x_t = x_t - \bar{x}_t$ and $\delta u_t = u_t - \bar{u}_t$ are state and control deviations, respectively, and $\delta Q_t(x_t, u_t) = Q_t(x_t, u_t) - Q_t(\bar{x}_t, \bar{u}_t)$ and $\delta V_t(x_t) = V_t(x_t) - V_t(\bar{x}_t)$.
The derivatives of $Q_t$ with respect to the nominal variables follow from~\eqref{eq:q-function} by the chain rule,
\begin{align}
    Q_{x,t} &= c_{x,t} + f_{x}^{\top}V_{x,t+1},\quad Q_{xx,t} = c_{xx,t} + f_{x}^{\top}V_{xx,t+1}f_{x},\\
    Q_{u,t} &= c_{u,t} + f_{u}^{\top}V_{x,t+1},\quad  Q_{uu,t} = c_{uu,t} + f_{u}^{\top}V_{xx,t+1}f_{u},\\[-1em]
    Q_{xu,t} &= c_{xu,t} + f_{x}^{\top}V_{xx,t+1}f_{u},
\end{align}
where we omit function arguments for a compact notation. The subscripts denote partial derivatives evaluated at the nominal variables, and second-order derivatives of the dynamics are omitted as in iLQR~\cite{tassa2012synthesis}. The minimization in~\eqref{eq:bellman} is over the mirror control, so the expansion has to be taken with respect to $y_t$. We denote derivatives of $Q_t(x_t,\nabla\psi^*(y_t))$ with respect to $y_t$ by the subscript $y$, for example $Q_{y,t}$ and $Q_{yy,t}$. Applying Lemma~\ref{lem:mirror-derivatives} to $Q_t(\bar{x}_t,\cdot)$ then gives the expansion with respect to $y_t$ in terms of the control-space terms. In the following, we write $\Psi_t \coloneqq \Psi(\bar{u}_t)$.

\begin{proposition}[Mirror-space expansion]
\label{prop:mirror-q-expansion}
The expansion of $Q_t(x_t,\nabla\psi^*(y_t))$ with respect to $y_t$ leads to
\begin{align}
    Q_{y,t} &= \Psi_t^{-1}\,Q_{u,t},\label{eq:mirror-q-gradient}\\
    Q_{xy,t} &= Q_{xu,t}\,\Psi_t^{-1},\label{eq:mirror-q-cross}\\
    Q_{yy,t} &= \Psi_t^{-1}\big[Q_{uu,t} - \mathcal{T}_\psi(\bar{u}_t)\cdot Q_{y,t}\big]\Psi_t^{-1},\label{eq:mirror-q-hessian}
\end{align}
while $Q_{x,t}$ and $Q_{xx,t}$ are unaffected.
\end{proposition}

\begin{proof}
Equations~\eqref{eq:mirror-q-gradient} and~\eqref{eq:mirror-q-hessian} are~\eqref{eq:mirror-gradient} and~\eqref{eq:mirror-hessian} with $g = Q_t(\bar{x}_t,\cdot)$. Applying~\eqref{eq:mirror-gradient} component-wise to $Q_{x,t}$ gives~\eqref{eq:mirror-q-cross}. The derivatives with respect to $x_t$ do not change.
\end{proof}

The third-order term in~\eqref{eq:mirror-q-hessian} is linear in $Q_{y,t}$ and therefore vanishes whenever $Q_{y,t}=0$. Following the iLQR Gauss--Newton approximation, we omit this term and use
\begin{equation}
    \label{eq:mirror-q-hessian-gn}
    Q_{yy,t} = \Psi_t^{-1}\,Q_{uu,t}\,\Psi_t^{-1},
\end{equation}
which is positive semidefinite whenever $Q_{uu,t}$ is, and which does not require $\mathcal{T}_\psi$. Section~\ref{sec:analysis} shows that, for linear-quadratic problems, this choice makes the backward pass compute a Levenberg--Marquardt step.

We apply Levenberg--Marquardt regularization, $\widehat{Q}_{yy,t} = Q_{yy,t} + \lambda I_{m}$, where $\lambda \geq 0$ is chosen such that $\widehat{Q}_{yy,t} \succ 0$ and is adapted with the gain-ratio rule of~\cite[Eq.~(2.21)]{madsen2004methods}. Minimizing the quadratic expansion of the $Q$-function with respect to the mirror control deviation $\delta y_t = y_t - \bar{y}_t$, for a given state deviation $\delta x_t$, yields the mirror feedback control law $\delta y_t^\star = k^y_t + K^y_t\,\delta x_t$, where $k^y_t = -\widehat{Q}_{yy,t}^{-1}\,Q_{y,t}, \, K^y_t = -\widehat{Q}_{yy,t}^{-1}\,Q_{xy,t}^\top$. Starting from $V_{x,T} = \nabla c_T(\bar{x}_T)$ and $V_{xx,T} = \nabla^2 c_T(\bar{x}_T)$, we propagate the value function backward as
\begin{equation}
    \label{eq:value-recursion}
    \begin{aligned}
        &V_{xx,t}
            = Q_{xx,t}
                + (K^y_t)^\top \widehat{Q}_{yy,t}K^y_t
             + (K^y_t)^\top Q_{xy,t}^\top
                + Q_{xy,t}K^y_t, \\
        &V_{x,t}
            = Q_{x,t}
                + Q_{xy,t}k^y_t
            + (K^y_t)^\top Q_{y,t}
                + (K^y_t)^\top \widehat{Q}_{yy,t}k^y_t, \\
        &\Delta V_t
            = Q_{y,t}^\top k^y_t
                + \frac{1}{2}\,(k^y_t)^\top \widehat{Q}_{yy,t}k^y_t.
    \end{aligned}
\end{equation}
Using $\widehat{Q}_{yy,t}$ throughout, the backward pass exactly minimizes the quadratic expansion augmented by the penalty $\frac{\lambda}{2}\|\delta y_t\|^2$ at every stage, as in~\cite{tassa2012synthesis}. The sum of $\Delta V_t$ over $t$ is the predicted change in cost, which enters the gain ratio used to adapt $\lambda$. The nominal trajectories are then updated through a feedback-controlled forward rollout. Starting from the initial state $x_0=x_s$, for $t=0,\ldots,T-1$, we compute
\begin{equation}
    \begin{split}
    \delta x_t &= x_t-\bar{x}_t, \quad
    y_t = \bar{y}_t+k^y_t+K^y_t\,\delta x_t,\\
    u_t &= \nabla\psi^*(y_t), \quad
    x_{t+1}=f(x_t,u_t).
    \end{split}\label{eq:forward-pass}
\end{equation}
The next proposition links this iteration to control-space iLQR.
\begin{proposition}[Relation to iLQR]
\label{prop:ilqr-relation}
Let $Q_{yy,t}$ be given by the Gauss--Newton approximation~\eqref{eq:mirror-q-hessian-gn}, and define $\widehat{Q}_{uu,t} \coloneqq Q_{uu,t} + \lambda\Psi_t^2$ together with the control-space gains $k^u_t = -\widehat{Q}_{uu,t}^{-1}\,Q_{u,t}, K^u_t = -\widehat{Q}_{uu,t}^{-1}\,Q_{xu,t}^\top$. Then, $k^y_t = \Psi_t\,k^u_t$ and $K^y_t = \Psi_t\,K^u_t$, and the value recursion~\eqref{eq:value-recursion} coincides with the control-space recursion of iLQR~\cite{tassa2012synthesis}, obtained from~\eqref{eq:value-recursion} by replacing the mirror-space terms with their control-space counterparts. Moreover, the forward rollout~\eqref{eq:forward-pass} applies
\begin{equation}
    \label{eq:mirror-update}
    u_t = \nabla\psi^*\big(\nabla\psi(\bar{u}_t) + \Psi_t\,\delta u_t\big),
\end{equation}
where $\delta u_t = k^u_t + K^u_t\,\delta x_t$ is the control-space step.
\end{proposition}

\begin{proof}
Under the Gauss--Newton approximation, the regularized mirror-space Hessian factors as $\widehat{Q}_{yy,t} = \Psi_t^{-1}\big(Q_{uu,t} + \lambda\Psi_t^2\big)\Psi_t^{-1} = \Psi_t^{-1}\,\widehat{Q}_{uu,t}\,\Psi_t^{-1}$, so that $\widehat{Q}_{yy,t}^{-1} = \Psi_t\,\widehat{Q}_{uu,t}^{-1}\,\Psi_t$. Inserting this factorization and the mirror-space derivatives of Proposition~\ref{prop:mirror-q-expansion} into the mirror gains yields $k^y_t = \Psi_t\,k^u_t,$
and $K^y_t = \Psi_t\,K^u_t$. In the value recursion, each mirror-space term then carries matching factors of $\Psi_t$ and $\Psi_t^{-1}$ that cancel, for example $(K^y_t)^\top \widehat{Q}_{yy,t}\,K^y_t = (K^u_t)^\top\Psi_t\,\Psi_t^{-1}\,\widehat{Q}_{uu,t}\,\Psi_t^{-1}\,\Psi_t\,K^u_t = (K^u_t)^\top \widehat{Q}_{uu,t}\,K^u_t,$
and the remaining terms cancel in the same way, which gives the control-space recursion. Finally, the forward rollout applies $y_t = \bar{y}_t + k^y_t + K^y_t\,\delta x_t = \nabla\psi(\bar{u}_t) + \Psi_t\,\delta u_t$.
\end{proof}

Proposition~\ref{prop:ilqr-relation} shows that the proposed method differs from iLQR in control space in two respects. First, the regularization $\lambda I_m$ is replaced by $\lambda\Psi_t^2$, so that step lengths are measured in mirror coordinates, since $\|\delta y_t\|^2 = \delta u_t^\top\Psi_t^2\,\delta u_t$ to first order. For the maps in Table~\ref{tab:mm-table}, $\Psi(u)$ grows without bound as $u$ approaches the boundary of $\mathcal{U}$, so that steps toward the boundary are damped increasingly strongly. Second, the control-space step is applied in mirror coordinates and mapped back by $\nabla\psi^*$. Since $\nabla\psi(\bar{u}_t) + \Psi_t\,\delta u_t$ is the first-order expansion of $\nabla\psi(\bar{u}_t+\delta u_t)$, the update~\eqref{eq:mirror-update} agrees with $u_t = \bar{u}_t + \delta u_t$ to first order, while $u_t\in\interior(\mathcal{U})$ holds for steps of any length.

\section{Convergence Analysis for Box-Constrained Linear-Quadratic Problems}\label{sec:analysis}
In this section, we analyze the convergence of the proposed method for box-constrained LQ control problems, where the mirror map introduces all nonlinearity into an otherwise convex quadratic program. Despite this nonlinearity, we show that the reparameterized problem is regular enough for a standard global convergence result of the Levenberg--Marquardt method to apply. In Section~\ref{sec:lin-quad-box}, we formulate the box-constrained problem in mirror coordinates. In Section~\ref{sec:lq-convergence}, we condense it into a nonlinear least-squares problem, relate the iLQR iteration to the Levenberg--Marquardt method, and verify the assumptions of the convergence result.

\subsection{Problem Formulation in Mirror Coordinates}
\label{sec:lin-quad-box}
We consider the linear-quadratic (LQ) objective
\begin{equation}
    \mathcal{J}_{\mathrm{LQ}}(u_{0:T-1}) = \frac{1}{2} x_T^\top P_T x_T + \frac{1}{2} \sum_{t=0}^{T-1} \Big[ x_t^\top P_t x_t + u_t^\top R_t u_t \Big],
\end{equation}
subject to dynamics $x_{t+1} = A x_t + B u_t$, initial state $x_0 = x_s$, and control box constraints $a \leq u_t \leq b$ for $t=0,\ldots,T-1$. We assume that $P_t\succeq 0$ and $R_t\succ 0$. We use the entropy DGF of Table~\ref{tab:mm-table}, defined on the interior of the box as
\begin{equation}
    \psi(u) = \sum_{i=1}^{m}\big[(b^i-u^i)\log(b^i-u^i) + (u^i-a^i)\log(u^i-a^i)\big].
\end{equation}
Its mirror map and inverse are
\begin{equation}
    y_t^i = \log\frac{u_t^i-a^i}{b^i-u_t^i}, \quad
    u_t^i = \left[\nabla \psi^*(y_t)\right]^i = \mu^i + \rho^i\tanh\!\left( y_t^i/2\right),
\end{equation}
where $\rho^i=(b^i-a^i)/2$ and $\mu^i= (b^i+a^i)/2$. The Hessian $\Psi(u_t)=\nabla^2\psi(u_t)$ is the diagonal, positive-definite matrix with entries  $\Psi^{ii} = 2\rho^i/((u_t^i-a^i)(b^i-u_t^i))$. Using the inverse map, we obtain an equivalent unconstrained reparameterization of the problem over the open box $\{u\colon a<u<b\}$. The cost in the new coordinates is
\begin{equation}
    \begin{split}
        &\mathcal{J}_{\mathrm{LQ}}([\nabla \psi^*(y_t)]_{t=0:T-1})= \widehat{\mathcal{J}}_{\mathrm{LQ}}(y_{0:T-1})\\
        &=\frac{1}{2} x_T^\top P_T x_T + \frac{1}{2}\!\sum_{t=0}^{T-1}\Big[x_t^\top P_t x_t + \nabla \psi^*(y_t)^\top R_t \nabla \psi^*(y_t) \Big].
    \end{split}
\end{equation}
The box-constrained LQ problem in mirror coordinates is then
\begin{equation}
    \begin{split}
        \minimize_{y_{0:T-1}} \quad &\widehat{\mathcal{J}}_{\mathrm{LQ}}(y_{0:T-1})\\
        \st \quad &x_{t+1} = A x_t + B\nabla \psi^*(y_t), \quad x_0=x_s.
    \end{split}\label{eq:mirror-lqr}
\end{equation}
The only remaining constraints are the dynamics, which we eliminate in the next section by condensing the problem.

\subsection{Convergence of the Levenberg--Marquardt Iteration}
\label{sec:lq-convergence}
The iLQR iteration with a regularized Hessian is closely related to the Levenberg--Marquardt method~\cite{nocedal2006numerical} applied to a nonlinear least-squares problem~\cite{tassa2012synthesis}. Hence, our analysis proceeds in three steps. First, we condense problem~\eqref{eq:mirror-lqr} into a nonlinear least-squares problem. Second, we specify the relation between iLQR and the Levenberg--Marquardt method. Third, we verify that the least-squares problem satisfies the assumptions of a global convergence result for trust-region methods~\cite[Thrm.~4.6]{nocedal2006numerical}.

For the first step, let us introduce the stacked state, control, and mirror-space variables, respectively
$U = \begin{bmatrix}u_0^\top & \cdots & u_{T-1}^\top\end{bmatrix}^\top$,
$X = \begin{bmatrix}x_1^\top & \cdots & x_{T}^\top\end{bmatrix}^\top$, and
$Y = \begin{bmatrix}y_0^\top & \cdots & y_{T-1}^\top\end{bmatrix}^\top$.
Throughout, $\|\cdot\|$ denotes the Euclidean norm of a vector and the induced spectral norm of a matrix. Furthermore, we denote the element-wise evaluation of the inverse map $\nabla \psi^*(\cdot)$ on the stacked mirror controls as $U(Y) = \begin{bmatrix}\nabla \psi^*(y_0)^\top & \cdots & \nabla \psi^*(y_{T-1})^\top\end{bmatrix}^\top.$ We can rewrite~\eqref{eq:mirror-lqr} in the form of a least-squares problem
\begin{equation}
\min_{Y}\,\phi(Y) \coloneqq \frac{1}{2}\,\|r(Y)\|^2,\label{eq:least-squares}
\end{equation}
where the residual is defined as $r(Y) = \Gamma U(Y) - \beta$ with
\begin{equation}
    \Gamma = \begin{bmatrix} P^{1/2} F \\ R^{1/2} \end{bmatrix}, \quad
    \beta = \begin{bmatrix} -P^{1/2} C\, x_s \\ 0 \end{bmatrix}.
\end{equation}
Here, $X = FU + Cx_s$ results from unrolling the dynamics, where $F$ is block lower-triangular with blocks $A^{t-1-\tau}B$ and $C$ stacks the powers $A,\ldots,A^T$. Moreover, $P = \blkdiag(P_1,\ldots,P_T)$ and $R = \blkdiag(R_0,\ldots,R_{T-1})$.

The residual is affine in $U(Y)$, so that all nonlinearity of~\eqref{eq:least-squares} stems from the stagewise inverse map. In particular, the Jacobian of the residual is $J_r(Y) =\Gamma J_U(Y)$, where
\begin{equation}
        J_U(Y) = \blkdiag\left(\nabla^2\psi^*(y_0),\dots,\nabla^2\psi^*(y_{T-1})\right),
\end{equation}
and $\nabla^2\psi^*(y_t) = \Psi(u_t)^{-1}$ with $u_t = \nabla\psi^*(y_t)$ by Proposition~\ref{prop:global-diffeo}, so that $J_U(Y)$ is diagonal.

For the second step, recall that Levenberg--Marquardt minimizes at each iterate $Y^{(k)}$ the Gauss--Newton model $\frac{1}{2}\|r(Y^{(k)}) + J_r(Y^{(k)})p\|^2$ within a trust region, solving $\big(J_r(Y^{(k)})^\top J_r(Y^{(k)}) + \lambda^{(k)} I\big)p^{(k)} = -J_r(Y^{(k)})^\top r(Y^{(k)})$ for some $\lambda^{(k)}\geq 0$~\cite{nocedal2006numerical}. This model ignores the second derivative of $U(Y)$, the third-order term omitted in~\eqref{eq:mirror-q-hessian-gn}. Thus, the backward pass of Section~\ref{sec:mirror-ilqr} with $\lambda=\lambda^{(k)}$ computes $p^{(k)}$, equivalent by Proposition~\ref{prop:ilqr-relation} to control-space iLQR with damping $\lambda\Psi_t^2$. The iteration differs only in applying feedback through $\nabla\psi^*$ during rollout and adapting $\lambda$ directly rather than via a trust-region radius.

For the third step, recall that~\cite[Thrm.~4.6]{nocedal2006numerical} requires that the objective is bounded from below, that its gradient is Lipschitz continuous, and that the model Hessians $J_r(Y^{(k)})^\top J_r(Y^{(k)})$ are uniformly bounded. The first property holds since $\phi(Y)\geq 0$ for all $Y\in\mathbb{R}^{Tm}$. Since $\nabla\phi = J_r^\top r$, the second property follows once $J_r$ and $r$ are bounded and Lipschitz continuous, and the third once $J_r$ is bounded. As $J_r(Y) = \Gamma J_U(Y)$ with a constant matrix $\Gamma$, the bounds on $J_r$ reduce to bounds on the first and second derivatives of the inverse map, which we establish in Lemma~\ref{lem:jacobian}. Lemma~\ref{lem:gradient-lipschitz} then provides the Lipschitz continuity of the gradient.

\begin{lemma}
    \label{lem:jacobian}
    For all $Y, Y_1, Y_2\in\mathbb{R}^{Tm}$, the Jacobian of the residual satisfies $\|J_r(Y)\| \leq \kappa$ and $\|J_r(Y_1) - J_r(Y_2)\| \leq L_r\,\|Y_1 - Y_2\|$, with $\kappa \coloneqq \|\Gamma\|\max_i \rho^i/2$ and $L_r \coloneqq \|\Gamma\|\max_i \rho^i/(3\sqrt{3})$.
\end{lemma}

\begin{proof}
    Since $J_r = \Gamma J_U$ with constant $\Gamma$, it suffices to bound $J_U$, whose norm is the maximum absolute value of its diagonal entries. The entry associated with stage $t$ and component $i$ is the derivative of $u^i(y) = \mu^i + \rho^i\tanh(y/2)$ evaluated at $[Y]_t^i$, where $\frac{\partial u^i}{\partial y} = \frac{\rho^i}{2}\sech^2(y/2), \frac{\partial^2 u^i}{\partial y^2} = -\frac{\rho^i}{2}\sech^2(y/2)\tanh(y/2)$. Since $0 < \sech^2(z) \leq 1$ for all $z\in\mathbb{R}$, the first derivative is bounded by $\rho^i/2$, hence $\|J_U(Y)\| \leq \max_i\rho^i/2$. The function $q(z) = \sech^2(z)\tanh(z)$ attains its maximum absolute value $2/(3\sqrt{3})$ where $q'(z) = \sech^2(z)\big(1 - 3\tanh^2(z)\big)$ vanishes, that is, at $\tanh^2(z) = 1/3$, so the second derivative is bounded by $\rho^i/(3\sqrt{3})$. By the mean value theorem, $\|J_U(Y_1) - J_U(Y_2)\| \leq \max_i\rho^i/(3\sqrt{3})\,\|Y_1 - Y_2\|$. Multiplying both bounds by $\|\Gamma\|$ completes the proof.
\end{proof}

\begin{lemma}
    \label{lem:gradient-lipschitz}
    The residual satisfies $\|r(Y)\| \leq \bar{r}$ for all $Y\in\mathbb{R}^{Tm}$, with $\bar{r} \coloneqq \|\Gamma\|\sup_{Y}\|U(Y)\| + \|\beta\| < \infty$, and the gradient $\nabla\phi(Y) = J_r(Y)^\top r(Y)$ is Lipschitz continuous with constant $L_r\,\bar{r} + \kappa^2$.
\end{lemma}

\begin{proof}
    Since each $u_t = \nabla\psi^*(y_t)$ lies in the box, $\sup_Y\|U(Y)\|$ is finite, and the bound on $r$ follows from the triangle inequality. By Lemma~\ref{lem:jacobian} and the mean value theorem, $r$ is Lipschitz continuous with constant $\kappa$. Hence, for any $Y_1, Y_2 \in \mathbb{R}^{Tm}$,
    \begin{align*}
            \|\nabla\phi(Y_1)-\nabla\phi(Y_2)\| &\leq \|J_r(Y_1)-J_r(Y_2)\| \|r(Y_1)\|\\
            &\quad + \|J_r(Y_2)\| \|r(Y_1)-r(Y_2)\|\\
            &\leq (L_r\bar{r}+\kappa^2) \|Y_1-Y_2\|.
    \end{align*}
\end{proof}

\begin{theorem}
    \label{thm:lq-convergence}
    Applied to~\eqref{eq:least-squares}, the Levenberg--Marquardt method in trust-region form~\cite[Alg.~4.1]{nocedal2006numerical} with a positive acceptance threshold generates iterates that satisfy $\lim_{k\to\infty}\nabla\phi(Y^{(k)}) = 0$.
\end{theorem}

\begin{proof}
    The claim follows from~\cite[Thrm.~4.6]{nocedal2006numerical}, whose assumptions we verify. The objective is bounded from below since $\phi\geq 0$, its gradient is Lipschitz continuous by Lemma~\ref{lem:gradient-lipschitz}, and the model Hessians satisfy $\|J_r(Y^{(k)})^\top J_r(Y^{(k)})\|\leq\kappa^2$ by Lemma~\ref{lem:jacobian}. Moreover, the trust-region step satisfies $\|p^{(k)}\|\leq\Delta^{(k)}$ and decreases the model at least as much as its minimizer along the negative gradient within the trust region.
\end{proof}

Theorem~\ref{thm:lq-convergence} guarantees that the iterates become stationary in mirror coordinates, regardless of whether constraints are active at the solution. Since $\nabla\phi(Y) = J_U(Y)\,\Gamma^\top\big(\Gamma U(Y)-\beta\big)$, where $\Gamma^\top(\Gamma U - \beta)$ is the gradient of the LQ objective in the stacked controls and $J_U(Y)$ is nonsingular for finite $Y$, this is stationarity in the original controls whenever the iterates remain bounded.

\begin{figure}[t]
    \centering
    \input{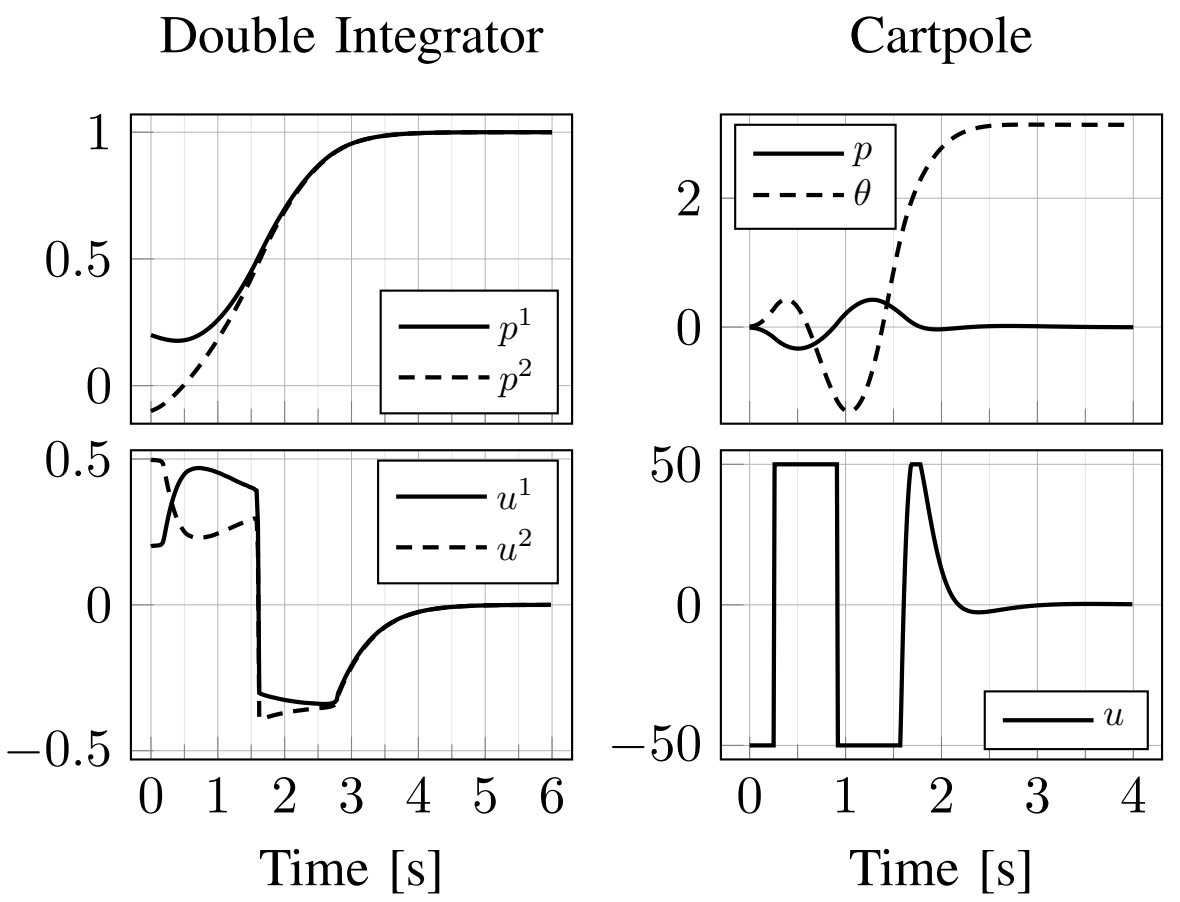}
    \vspace{-0.25cm}
    \caption{MPC closed-loop trajectories. Left: double integrator under
    polytopic constraints. Right: cart-pole swing-up under box constraints.
    Top: states; bottom: controls.}
    \label{fig:mpc-sim}
\end{figure}

\begin{figure}[t]
    \centering
    \input{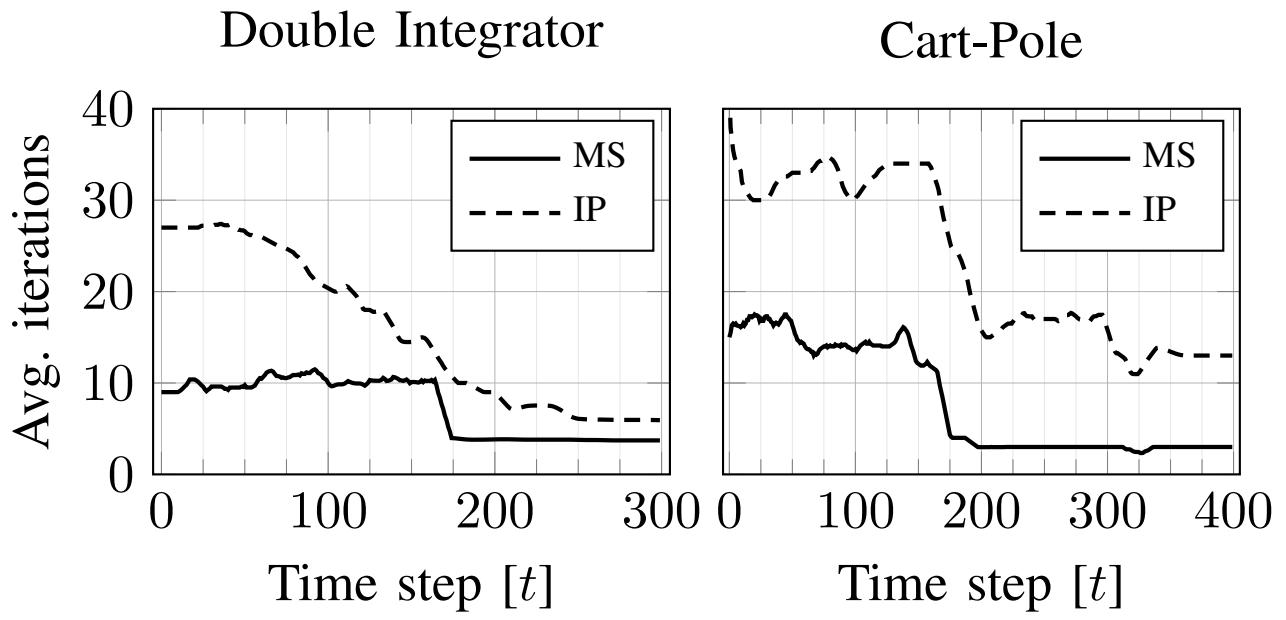}
    \vspace{-0.25cm}
    \caption{Number of iterations per MPC step, averaged over 10 initial conditions, for the proposed mirror-space (MS) method and interior-point (IP) differential dynamic programming. The left and right panels show the double-integrator and cart-pole tasks, respectively.}
    \label{fig:iterations}
\end{figure}

\section{Experiments}
\label{sec:experiments}
In this section, we evaluate the proposed method on a linear and a nonlinear MPC task and compare its empirical convergence against interior-point differential dynamic programming (DDP)~\cite{pavlov2021interior}. The implementation details are available in the supplementary code repository\footnote{\href{https://github.com/casiacob/mirror-space-ilqr}{https://github.com/casiacob/mirror-space-ilqr}}.

For our evaluation, we consider the following two systems:
\paragraph{\textup{\textbf{Double integrator with polytopic constraints}}}
For the linear case, we consider the double integrator in two dimensions~\cite[Ex.~4.2]{sarkka2023bayesian}. The state $x=\begin{bmatrix} p^1 & p^2 & v^1 & v^2 \end{bmatrix}^\top$ includes the positions and velocities, and the control $u=\begin{bmatrix} u^1 & u^2 \end{bmatrix}^\top$ comprises the accelerations. The objective is to drive the system from the initial state $x_0=\begin{bmatrix} 0.2&-0.1&-0.1&0.1 \end{bmatrix}^\top$ to the desired state $x_d=\begin{bmatrix} 1&1&0&0 \end{bmatrix}^\top$. The controls are subject to the polytopic constraints $|u^1| \leq 0.5$, $|u^2| \leq 0.5$, $|u^1+u^2| \leq 0.7$, and $|u^1-u^2| \leq 0.7$, which define the convex set $\mathcal{U}=\{u\colon Gu \leq h\}$ with $G\in\mathbb{R}^{8\times 2}$ and $h\in\mathbb{R}^{8}$. We choose the DGF $\psi(u) = \frac{1}{2} u^\top u - \sum_{j=1}^{8} \log s^j(u),$ where $s(u) = h - G u$ is the slack vector. The corresponding mirror map is $y = \nabla \psi (u) = u + G^\top \big[1/s^j(u)\big]_{j=1:8}$. The inverse map $\nabla\psi^*$ is intractable. Hence, we compute $u = \nabla\psi^*(y)$ by solving $\nabla\psi(u) = y$ with Newton's root-finding method~\cite{nocedal2006numerical}.

\paragraph{\textup{\textbf{Cart-pole with box constraints}}}
For the nonlinear case, we consider the cart-pole system~\cite{spong2005underactuated}. The state vector includes the cart position, pole angle, and their velocities $x=\begin{bmatrix} p&\theta& v&\omega \end{bmatrix}^\top$, respectively. The control vector is the force applied to the cart. We consider a swing-up task, where the objective is to drive the system from the initial state $x_0 = \begin{bmatrix} -0.01&0.01&-0.01&0.01 \end{bmatrix}^\top$ to the desired state $x_d = \begin{bmatrix} 0&\pi&0&0 \end{bmatrix}^\top$. The force is bounded such that $u\in[-50,50]$. We use the box (entropy) DGF from Table~\ref{tab:mm-table}.

Fig.~\ref{fig:mpc-sim} shows the MPC state and control trajectories for the double integrator and cart-pole. Both systems reach their targets while satisfying the control constraints.

Fig.~\ref{fig:iterations} compares the proposed mirror-space (MS) method and interior point DDP (IP) iterations per MPC step, averaged over 10 random initial conditions per task. Each iteration comprises one backward and one forward pass. MS requires substantially fewer iterations on both tasks.

\section{Conclusion}
We developed a dynamic programming method for nonlinear optimal control with convex control constraints using a Legendre-type DGF to transform constrained controls into unconstrained mirror coordinates. Running iLQR in these coordinates requires the Gauss--Newton model of the Bellman recursion in the mirror controls, and hence their derivatives. Using convex conjugacy, we express these derivatives in terms of control-space derivatives and the Hessian of the DGF, without differentiating the inverse map. The resulting algorithm is a damped iLQR in mirror coordinates, which keeps every iterate strictly feasible. For box-constrained linear-quadratic problems, we proved global convergence of the Levenberg--Marquardt iteration. In MPC experiments, our method requires substantially fewer iterations than interior-point methods.

\section*{Acknowledgement}
We used generative AI tools to assist with mathematical proofs, method implementation, and manuscript editing. We did not use them for conceptual development, hypothesis formulation, research design, interpretation of results, or literature discovery and summarization. No sections were explicitly written by AI tools. We reviewed all AI-assisted work and take full responsibility for the final content of this paper.

\begin{appendices}
\section{Proof of Lemma~\ref{lem:mirror-derivatives}}
\label{app:proof-mirror-derivatives}
Let $W(y) \coloneqq \Psi(\nabla\psi^*(y))^{-1}$, so that $\partial u/\partial y = W$ by Proposition~\ref{prop:global-diffeo}. The chain rule gives $\nabla h(y) = (\partial u/\partial y)^\top\nabla g(u) = W\nabla g(u)$, which is~\eqref{eq:mirror-gradient}. Differentiating the $i$-th component of $W\nabla g(u)$ with respect to $y^\ell$ by the product rule yields
\begin{equation}
    [\nabla^2 h]_{i\ell} = \underbrace{\sum_j \frac{\partial W_{ij}}{\partial y^\ell}\,[\nabla g]^j}_{[E_1]_{i\ell}} + \underbrace{\sum_j W_{ij}\,\frac{\partial [\nabla g]^j}{\partial y^\ell}}_{[E_2]_{i\ell}}.
\end{equation}
For $E_2$, the chain rule gives $\partial[\nabla g]^j/\partial y^\ell = \sum_k [\nabla^2 g]_{jk}\,W_{k\ell}$, hence $E_2 = W\,\nabla^2 g(u)\,W$. For $E_1$, we apply the identity $d(H^{-1}) = -H^{-1}(dH)H^{-1}$ to $H = \Psi(u)$ and use $\partial\Psi_{pq}/\partial u^k = [\mathcal{T}_\psi]_{pqk}$ from~\eqref{eq:third-order-tensor}, which gives
\begin{equation}
    \frac{\partial W_{ij}}{\partial u^k} = -\sum_{p,q} W_{ip}\,[\mathcal{T}_\psi]_{pqk}\,W_{qj}.
\end{equation}
Together with $\partial u^k/\partial y^\ell = W_{k\ell}$, the chain rule yields
\begin{equation}
    \frac{\partial W_{ij}}{\partial y^\ell} = -\sum_{k,p,q} W_{ip}\,[\mathcal{T}_\psi]_{pqk}\,W_{qj}\,W_{k\ell}.
\end{equation}
Multiplying by $[\nabla g]^j$ and summing over $j$, we recognize $\sum_j W_{qj}[\nabla g]^j = [\nabla h]^q$ from~\eqref{eq:mirror-gradient}, so that
\begin{equation}
    \begin{split}
        [E_1]_{i\ell} &= -\sum_{k,p} W_{ip}\Big(\sum_q [\mathcal{T}_\psi]_{pqk}\,[\nabla h]^q\Big) W_{k\ell}\\
        &= -\big[W\,\big(\mathcal{T}_\psi\cdot\nabla h\big)\,W\big]_{i\ell},
    \end{split}
\end{equation}
where the inner sum equals~\eqref{eq:tensor-dot} by symmetry of $\mathcal{T}_\psi$ for $\psi\in C^3$. Adding $E_1$ and $E_2$ yields~\eqref{eq:mirror-hessian}.
\end{appendices}

\bibliographystyle{IEEEtran}
\bibliography{Ref}

\end{document}